\documentclass{article}
\usepackage{amssymb}
\usepackage{amsfonts}
\usepackage{amsmath}

\begin{document}

\begin{center}
{\Large Some remarks regarding }$(a,b,x_{0},x_{1})-${\Large numbers \ and} $%
(a,b,x_{0},x_{1})-${\Large quaternions}

\begin{equation*}
\end{equation*}%
Cristina FLAUT and Diana SAVIN 
\begin{equation*}
\end{equation*}
\end{center}

\textbf{Abstract. }{\small In this paper we define and study properties and
applications of }$(a,b,x_{0},x_{1})-${\small elements in some special cases. 
}%
\begin{equation*}
\end{equation*}


\medskip

\begin{equation*}
\end{equation*}

\bigskip

\textbf{1. Introduction}%
\begin{equation*}
\end{equation*}%
\newline
\qquad We consider the following difference equation of degree two%
\begin{equation}
d_{n}=ad_{n-1}+bd_{n-2},d_{0}=x_{0},d_{1}=x_{1}.  \tag{1.1}
\end{equation}

The study of this equation in the general case could be interesting from a
pure mathematical point of view, but some particular cases produced very
good and important applications. This is the reason for which the properties
and applications of this difference equation were intensively studied in the
last years, in some important special cases.

We will define the $(a,b,x_{0},x_{1})-$\textit{numbers} to be the numbers
which satisfy the equations $\left( 1.1\right) ,$ where $a,b,x_{0},x_{1}$
are arbitrary integer numbers. For example, if we consider $a=b=1,$ $%
x_{0}=0,x_{1}=1,$ we obtain the Fibonacci numbers and if we consider $a=b=1,$
$x_{0}=2,x_{1}=1,$ we obtain the Lucas numbers. The properties and
applications of some particulae cases of these numbers are various and were
extended to other algebraic structures, or used as applications in the
Coding Theory. \ We refer here to a small part of those papers which
approached this subject: [Ak, Ko, To; 14], [Ca; 15], [Fa, Pl; 07(1)], [Fa,
Pl; 07(2)], [Fa, Pl; 09], [Fl, Sa; 15], [Fl, Sh;13], [Fl, Sh;15], [Fl, Sh,
Vl; 17], [ Sa; 17], [Gu, Nu; 15], [Ha; 12], [Ho; 63], [Na, Ha; 09], [Ra;
15], [St; 06], [Sw; 73].

In this paper, we will provide some properties and applications of the%
\textbf{\newline
}$(a,b,x_{0},x_{1})-$elements, in some special cases.

This paper is organized as follow: in Section 3, we found the generating
function for the $\left( 1,1,p+2q,q\right) -$numbers (the generalized
Fibonacci-Lucas numbers), the square of this function, Cassini identity for
these elements and other interesting properties and applications. In Section
4, we found the generating function for the $\left( 1,1,p+2q,q\right) -$%
quaternions (generalized Fibonacci-Lucas quaternions), Binet's formula,
Catalan's and Cassini's identities. In Section 5, we studied the $\left(
1,a,0,1\right) -$numbers, the $\left( 1,a,2,1\right) -$numbers, the\newline
$\left( 1,a,p+2q,q\right) -$quaternions and we gave some interesting
properties, as for example an algebraic structure for the last of them.

\begin{equation*}
\end{equation*}

\bigskip

\textbf{2. Preliminaries }%
\begin{equation*}
\end{equation*}%
\qquad \qquad

First of all, we recall some elementary properties of the Fibonacci and
Lucas numbers, properties which will be necessary in the proofs of this
paper.\newline
Let $(f_{n})_{n\geq 0}$ be the Fibonacci sequence and let $(l_{n})_{n\geq 0}$
be the Lucas sequence. Let $\alpha =\frac{1+\sqrt{5}}{2}$ and $\beta =\frac{%
1-\sqrt{5}}{2}.$\newline
\smallskip \newline
\textbf{Binet's formula for Fibonacci sequence.} 
\begin{equation*}
f_{n}=\frac{\alpha ^{n}-\beta ^{n}}{\alpha -\beta }=\frac{\alpha ^{n}-\beta
^{n}}{\sqrt{5}},\ \ \left( \forall \right) n\in \mathbb{N}.
\end{equation*}%
\textbf{Binet's formula for Lucas sequence.} 
\begin{equation*}
l_{n}=\alpha ^{n}+\beta ^{n},\ \ \left( \forall \right) n\in \mathbb{N}.
\end{equation*}

\medskip

\textbf{Proposition 2.1.} ([Fib.]). \textit{Let} $(f_{n})_{n\geq 0}$ \textit{%
be the Fibonacci sequence} \textit{and let } $(l_{n})_{n\geq 0}$ \textit{be
the Lucas sequence.} \textit{The following properties hold:}\newline
i) 
\begin{equation*}
f_{n}^{2}+f_{n+1}^{2}=f_{2n+1},\forall ~n\in \mathbb{N};
\end{equation*}%
ii) 
\begin{equation*}
f_{n+1}^{2}-f_{n-1}^{2}=f_{2n},\forall ~n\in \mathbb{N}^{\ast };
\end{equation*}%
iii) 
\begin{equation*}
l_{n}^{2}-f_{n}^{2}=4f_{n-1}f_{n+1},\forall ~n\in \mathbb{N}^{\ast };
\end{equation*}%
iv) 
\begin{equation*}
l_{n}^{2}+l_{n+1}^{2}=5f_{2n+1},\forall ~n\in \mathbb{N};
\end{equation*}%
v) 
\begin{equation*}
l_{n}^{2}=l_{2n}+2\left( -1\right) ^{n},\forall ~n\in \mathbb{N}^{\ast };
\end{equation*}%
vi) 
\begin{equation*}
f_{n+1}+f_{n-1}=l_{n},\forall ~n\in \mathbb{N}^{\ast };
\end{equation*}%
vii) 
\begin{equation*}
l_{n}+l_{n+2}=5f_{n+1},\forall ~n\in \mathbb{N};
\end{equation*}%
viii) 
\begin{equation*}
f_{n}+f_{n+4}=3f_{n+2},\forall ~n\in \mathbb{N};
\end{equation*}%
ix) 
\begin{equation*}
f_{m}l_{m+p}=f_{2m+p}+\left( -1\right) ^{m+1}f_{p},\forall ~m,p\in \mathbb{N}%
;
\end{equation*}%
x) 
\begin{equation*}
f_{m+p}l_{m}=f_{2m+p}+\left( -1\right) ^{m}f_{p},\forall ~m,p\in \mathbb{N};
\end{equation*}%
xi) 
\begin{equation*}
f_{m}f_{m+p}=\frac{1}{5}\left( l_{2m+p}+\left( -1\right) ^{m+1}l_{p}\right)
,\forall ~m,p\in \mathbb{N};
\end{equation*}%
xii) 
\begin{equation*}
l_{m}l_{p}+5f_{m}f_{p}=2l_{m+p},\forall ~m,p\in \mathbb{N}.
\end{equation*}%
\smallskip \newline
Let $\left( a_{n}\right) _{n\geq 1}$ and $\left( b_{n}\right) _{n\geq 1}$ be
two sequences of real numbers and let $A\left( z\right) =\sum\limits_{n\geq
1}a_{n}z^{n}$ and $B\left( z\right) =\sum\limits_{n\geq 1}b_{n}z^{n}$ be
their generating functions. We recall that their product has the form\newline
$A\left( z\right) B\left( z\right) =\sum\limits_{n\geq 1}s_{n}z^{n},$ where $%
s_{n}=\sum\limits_{k=1}^{n}a_{k}b_{n-k}.$ 
\begin{equation*}
\end{equation*}%
\textbf{\ 3. Generalized Fibonacci- Lucas numbers. Some properties and
applications}%
\begin{equation*}
\end{equation*}

Let $\left( f_{n}\right) _{n\geq 0}$ be the Fibonacci sequence 
\begin{equation*}
f_{n}=f_{n-1}+f_{n-2},\;n\geq 2,f_{0}=0;f_{1}=1,
\end{equation*}%
and $\left( l_{n}\right) _{n\geq 0}$ be the Lucas sequence\textbf{\ }%
\begin{equation*}
l_{n}=l_{n-1}+l_{n-2},\;n\geq 2,l_{0}=2;l_{1}=1.
\end{equation*}%
In the paper [Fl, Sa; 15], we introduced the generalized Fibonacci-Lucas
numbers. If $n$ is an arbitrary positive integer and $p,q$ are two arbitrary
integers, the sequence $\left( g_{n}\right) _{n\geq 1},$ where\textbf{\ }%
\begin{equation*}
g_{n+1}=pf_{n}+ql_{n+1},\;n\geq 0,
\end{equation*}
with $g_{0}=p+2q,~$is called the generalized Fibonacci-Lucas numbers. We
remark that $g_{n}=g_{n-1}+g_{n-2},$ with $g_{0}=p+2q,g_{1}=q$ and these
numbers are actually the $\left( 1,1,p+2q,q\right) -$numbers.

To avoid confusions, in the following, we will use the notation $g_{n}^{p,q}$
instead of $g_{n}.$

Let $\left( g_{n}^{p,q}\right) _{n\geq 1}$ be the generalized Fibonacci-
Lucas numbers and let $A$ be the generating function for these numbers 
\begin{equation*}
A\left( z\right) =\sum\limits_{n\geq 1}g_{n}^{p,q}z^{n}.
\end{equation*}
In the next proposition we determine this function.\newline
\medskip

\textbf{Proposition 3.1.} \textit{With the above notations, the following
relation is true:} 
\begin{equation*}
A\left( z\right) =\frac{qz+\left( p+2q\right) z^{2}}{1-z-z^{2}}.
\end{equation*}%
\textbf{Proof.} 
\begin{equation}
A\left( z\right) =g_{1}^{p,q}z+g_{2}^{p,q}z^{2}+...g_{n}^{p,q}z^{n}+...\ \  
\tag{3.1}
\end{equation}%
\begin{equation}
zA\left( z\right)
=g_{1}^{p,q}z^{2}+g_{2}^{p,q}z^{3}+...g_{n-1}^{p,q}z^{n}+...\ \   \tag{3.2}
\end{equation}%
\begin{equation}
z^{2}A\left( z\right)
=g_{1}^{p,q}z^{3}+g_{2}^{p,q}z^{4}+...g_{n-2}^{p,q}z^{n}+...\ \ .  \tag{3.3}
\end{equation}%
By adding the equalities (3.2) and (3.3) member by member, we obtain 
\begin{equation*}
A\left( z\right) \left( 1-z-z^{2}\right) =qz+\left( p+2q\right) z^{2},
\end{equation*}%
therefore%
\begin{equation*}
A\left( z\right) =\frac{qz+\left( p+2q\right) z^{2}}{1-z-z^{2}}.
\end{equation*}%
$\Box \smallskip $\medskip

\textbf{Proposition 3.2.} \textit{If} $A^{2}\left( z\right)
=\sum\limits_{n\geq 1}s_{n}z^{n},$ \textit{then} 
\begin{equation*}
5s_{n}=ng_{n}^{10pq,p^{2}+5q^{2}}+g_{n}^{p^{2}+5q^{2}-10pq,5pq}+g_{n-1}^{p^{2}+5q^{2},0}-ng_{n-1}^{0,p^{2}}.
\end{equation*}%
\textbf{Proof.} We know that $s_{n}=\sum%
\limits_{k=1}^{n}g_{k}^{p,q}g_{n-k}^{p,q}.$ This equality is equivalent with 
\begin{equation}
s_{n}=\sum\limits_{k=1}^{n}\left(
p^{2}f_{k-1}f_{n-k-1}+pqf_{k-1}l_{n-k}+pqf_{n-k-1}l_{k}+q^{2}l_{k}l_{n-k}%
\right) \ \   \tag{3.4}
\end{equation}%
Using Binet's formulas for Fibonacci and Lucas numbers, we have 
\begin{equation*}
pqf_{k-1}l_{n-k}+pqf_{n-k-1}l_{k}=
\end{equation*}%
\begin{equation*}
=\frac{pq}{\sqrt{5}}\left[ \left( \alpha ^{k-1}-\beta ^{k-1}\right) \left(
\alpha ^{n-k}+\beta ^{n-k}\right) +\left( \alpha ^{n-k-1}-\beta
^{n-k-1}\right) \left( \alpha ^{k}+\beta ^{k}\right) \right] =
\end{equation*}%
\begin{equation*}
=\frac{pq}{\sqrt{5}}\left[ 2\alpha ^{n-1}-2\beta ^{n-1}+\alpha ^{k-1}\beta
^{n-k-1}\left( \beta -\alpha \right) +\alpha ^{n-k-1}\beta ^{k-1}\left(
\beta -\alpha \right) \right] =
\end{equation*}%
\begin{equation*}
=pq\left( 2\frac{\alpha ^{n-1}-\beta ^{n-1}}{\sqrt{5}}-\alpha ^{k-1}\beta
^{n-k-1}-\alpha ^{n-k-1}\beta ^{k-1}\right) .
\end{equation*}%
Using again Binet's formula for Fibonacci numbers, we obtain 
\begin{equation}
pqf_{k-1}l_{n-k}+pqf_{n-k-1}l_{k}=pq\left( 2f_{n-1}-\alpha ^{k-1}\beta
^{n-k-1}-\alpha ^{n-k-1}\beta ^{k-1}\right) \   \tag{3.5}
\end{equation}%
From (3.5), it results that 
\begin{equation*}
\sum\limits_{k=1}^{n}\left( pqf_{k-1}l_{n-k}+pqf_{n-k-1}l_{k}\right) =
\end{equation*}%
\begin{equation*}
=2pqnf_{n-1}-pq\sum\limits_{k=1}^{n}\alpha ^{k-1}\beta
^{n-k-1}-pq\sum\limits_{k=1}^{n}\alpha ^{n-k-1}\beta ^{k-1}=
\end{equation*}%
\begin{equation*}
=2pqnf_{n-1}-pq\left( \frac{\alpha ^{n-1}-\beta ^{n-1}}{\alpha -\beta }%
+\alpha ^{n-1}\beta ^{-1}\right) -pq\left( \frac{\alpha ^{n-1}-\beta ^{n-1}}{%
\alpha -\beta }+\alpha ^{-1}\beta ^{n-1}\right) =
\end{equation*}%
\begin{equation*}
=2pqnf_{n-1}-pq\left( f_{n-1}+\alpha ^{n-1}\beta ^{-1}\right) -pq\left(
f_{n-1}+\alpha ^{-1}\beta ^{n-1}\right) =
\end{equation*}%
\begin{equation*}
=2pqnf_{n-1}-2pqf_{n-1}-pq\frac{\alpha ^{n}+\beta ^{n}}{\alpha \cdot \beta }%
=2pqnf_{n-1}-2pqf_{n-1}+pq\left( \alpha ^{n}+\beta ^{n}\right) =
\end{equation*}%
\begin{equation*}
=2pqf_{n-1}\left( n-1\right) +pql_{n}.
\end{equation*}%
Therefore, we get 
\begin{equation}
\sum\limits_{k=1}^{n}\left( pqf_{k-1}l_{n-k}+pqf_{n-k-1}l_{k}\right)
=2pqf_{n-1}\left( n-1\right) +pql_{n}\ \   \tag{3.6}
\end{equation}%
We use again Binet's formula for Fibonacci and Lucas numbers and we obtain 
\begin{equation*}
\sum\limits_{k=1}^{n}p^{2}f_{k-1}f_{n-k-1}=p^{2}\sum\limits_{k=1}^{n}\frac{%
\alpha ^{k-1}-\beta ^{k-1}}{\sqrt{5}}\frac{\alpha ^{n-k-1}-\beta ^{n-k-1}}{%
\sqrt{5}}=
\end{equation*}%
\begin{equation*}
=\frac{p^{2}}{5}\sum\limits_{k=1}^{n}\left( \alpha ^{n-2}+\beta
^{n-2}\right) -\frac{p^{2}}{5}\left[ \sum\limits_{k=1}^{n}\left( \alpha
^{k-1}\beta ^{n-k-1}+\alpha ^{n-k-1}\beta ^{k-1}\right) \right] =
\end{equation*}%
\begin{equation*}
=\frac{p^{2}}{5}nl_{n-2}-\frac{p^{2}}{5}\left[ \sum\limits_{k=1}^{n}\left(
\alpha ^{k-1}\beta ^{n-k-1}+\alpha ^{n-k-1}\beta ^{k-1}\right) \right] .
\end{equation*}%
From relation (3.6) we have 
\begin{equation*}
\sum\limits_{k=1}^{n}\left( \alpha ^{k-1}\beta ^{n-k-1}+\alpha ^{n-k-1}\beta
^{k-1}\right) =2f_{n-1}+\frac{\alpha ^{n}+\beta ^{n}}{\alpha \beta }%
=2f_{n-1}-l_{n}.
\end{equation*}%
Therefore we obtain the following relation 
\begin{equation*}
\sum\limits_{k=1}^{n}p^{2}f_{k-1}f_{n-k-1}=\frac{p^{2}n}{5}l_{n-2}-\frac{%
p^{2}}{5}\left( 2f_{n-1}-l_{n}\right) .
\end{equation*}%
Applying Proposition 2.1 (vi), we get 
\begin{equation}
\sum\limits_{k=1}^{n}p^{2}f_{k-1}f_{n-k-1}=\frac{p^{2}\left(
nl_{n-2}+f_{n}\right) }{5}.  \tag{3.7}
\end{equation}%
Using Binet's formula for Lucas numbers and after for Fibonacci numbers, we
have the following relation 
\begin{equation*}
\sum\limits_{k=1}^{n}q^{2}l_{k}l_{n-k}=q^{2}\sum\limits_{k=1}^{n}\left(
\alpha ^{k}+\beta ^{k}\right) \left( \alpha ^{n-k}+\beta ^{n-k}\right) =
\end{equation*}%
\begin{equation*}
=q^{2}\sum\limits_{k=1}^{n}\left( \alpha ^{n}+\beta ^{n}\right)
+q^{2}\sum\limits_{k=1}^{n}\left( \alpha ^{k}\beta ^{n-k}+\alpha ^{n-k}\beta
^{k}\right) =
\end{equation*}%
\begin{equation*}
=q^{2}nl_{n}+q^{2}\left[ -\left( \alpha ^{n}+\beta ^{n}\right) +2\frac{%
\alpha ^{n+1}-\beta ^{n+1}}{\alpha -\beta }\right] =
\end{equation*}%
\begin{equation*}
=q^{2}nl_{n}+q^{2}\left( -l_{n}+2f_{n+1}\right) .
\end{equation*}%
Using Proposition 2.1 (vi), we obtain 
\begin{equation}
\sum\limits_{k=1}^{n}q^{2}l_{k}l_{n-k}=q^{2}\left( nl_{n}+f_{n}\right) . 
\tag{3.8}
\end{equation}%
Adding member by member the equalities (3.6), (3.7), (3.8), we have 
\begin{equation*}
s_{n}=2pqf_{n-1}\left( n-1\right) +pql_{n}+\frac{p^{2}\left(
nl_{n-2}+f_{n}\right) }{5}+q^{2}\left( nl_{n}+f_{n}\right) =
\end{equation*}%
\begin{equation*}
=2pqf_{n-1}\left( n-1\right) +pql_{n}+\frac{p^{2}\left(
nl_{n}-nl_{n-1}+f_{n-2}+f_{n-1}\right) }{5}+q^{2}\left(
nl_{n}+f_{n-2}+f_{n-1}\right) =
\end{equation*}%
\begin{equation*}
=\left[ pq\left( 2n-2\right) \text{+}\left( \frac{p^{2}}{5}+q^{2}\right) %
\right] f_{n-1}\text{+}\left( \frac{np^{2}}{5}+pq+nq^{2}\right) l_{n}\text{+}%
\left( \frac{p^{2}}{5}+q^{2}\right) f_{n-2}-\frac{np^{2}}{5}l_{n-1}.
\end{equation*}%
It results that 
\begin{equation*}
5s_{n}=n\left[ 10pqf_{n-1}+\left( p^{2}+5q^{2}\right) l_{n}\right] +
\end{equation*}%
\begin{equation*}
+\left( p^{2}+5q^{2}-10pq\right) f_{n-1}+5pql_{n}+\left( p^{2}+5q^{2}\right)
f_{n-2}-np^{2}l_{n-1}.
\end{equation*}%
Finally, we obtain 
\begin{equation*}
5s_{n}=ng_{n}^{10pq,p^{2}+5q^{2}}+g_{n}^{p^{2}+5q^{2}-10pq,5pq}+g_{n-1}^{p^{2}+5q^{2},0}-ng_{n-1}^{0,p^{2}}.
\end{equation*}%
$\Box \smallskip $

For Fibonacci and Lucas numbers, there are very well known the Cassini's
identities, namely 
\begin{equation*}
f_{n+1}f_{n-1}-f_{n}^{2}=\left( -1\right) ^{n},\;\left( \forall \right)
\;n\in \mathbb{N}^{\ast },
\end{equation*}%
for Fibonacci numbers and 
\begin{equation*}
l_{n+1}l_{n-1}-l_{n}^{2}=5\cdot \left( -1\right) ^{n-1},\;\left( \forall
\right) \;n\in \mathbb{N}^{\ast },
\end{equation*}%
for Lucas numbers. In the following, we extend these results for the
generalized Fibonacci-Lucas numbers.\medskip

\textbf{Proposition 3.3.} \textit{Let} $p,q$ \textit{be arbitrary integers.
Then, we have:} 
\begin{equation*}
g_{n+1}^{p,q}g_{n-1}^{p,q}-\left( g_{n}^{p,q}\right) ^{2}=\left( -1\right)
^{n-1}\left[ p^{2}+5q^{2}+5pq\right] ,\text{ for all}\;n\in \mathbb{N}%
,\;n\geq 2.
\end{equation*}%
\textbf{Proof. }\ We have that 
\begin{equation*}
g_{n+1}^{p,q}g_{n-1}^{p,q}-\left( g_{n}^{p,q}\right) ^{2}=\left(
pf_{n}+ql_{n+1}\right) \left( pf_{n-2}+ql_{n-1}\right) -\left(
pf_{n-1}+ql_{n}\right) ^{2}=
\end{equation*}%
\begin{equation*}
=p^{2}\left( f_{n-2}f_{n}-f_{n-1}^{2}\right) +q^{2}\left(
l_{n-1}l_{n+1}-l_{n}^{2}\right) +pq\left(
f_{n-2}l_{n+1}+f_{n}l_{n-1}-2f_{n-1}l_{n}\right)
\end{equation*}%
Using Cassini identities for Fibonacci and Lucas numbers and Proposition 2.1
(ix; x), we obtain

\begin{equation*}
g_{n+1}^{p,q}g_{n-1}^{p,q}-\left( g_{n}^{p,q}\right) ^{2}=
\end{equation*}%
\begin{equation*}
=p^{2}\left( -1\right) ^{n-1}+5q^{2}\left( -1\right) ^{n-1}+
\end{equation*}%
\begin{equation*}
+pq\left[ f_{2n-1}+\left( -1\right) ^{n-1}f_{3}+f_{2n-1}+\left( -1\right)
^{n-1}f_{1}-2f_{2n-1}-2\left( -1\right) ^{n}f_{1}\right] =
\end{equation*}%
\begin{equation*}
=p^{2}\left( -1\right) ^{n-1}+5q^{2}\left( -1\right) ^{n-1}+5pq\left(
-1\right) ^{n-1}=
\end{equation*}%
\begin{equation*}
=\left( -1\right) ^{n-1}\left[ p^{2}+5q^{2}+5pq\right] .
\end{equation*}%
$\Box \smallskip $

For the generalized Fibonacci-Lucas numbers, we introduce the following
quadratic matrix 
\begin{equation}
M_{n}=\left( 
\begin{array}{cc}
g_{n+1}^{p,q} & g_{n}^{p,q} \\ 
g_{n}^{p,q} & g_{n-1}^{^{p,q}}%
\end{array}%
\right) .  \tag{3.9}
\end{equation}

\textbf{Proposition 3.4.} \textit{With the above notations,} \textit{the
following affirmations are true:}

\textit{i)} $M_{n}=M_{n-1}+M_{n-2}.$

\textit{ii) For} $n\in \mathbb{N},\;n\geq 2,$ \textit{we have} det$M=\left(
-1\right) ^{n-1}\left[ p^{2}+5q^{2}+5pq\right] .\medskip $

\textbf{Proof.} \newline
i) Since $g_{n+1}^{p,q}=g_{n}^{p,q}+g_{n-1}^{p,q}$, $%
g_{n}^{p,q}=g_{n-1}^{p,q}+g_{n-2}^{p,q},~g_{n-1}^{p,q}=g_{n-2}^{p,q}+g_{n-3}^{p,q}, 
$ we obtain the asked relation.\newline
ii) We use Proposition 3.3. $\Box \smallskip $

\textbf{Remark 3.5.} In [St; 06] and [Ba, Pr; 09], were presented some
applications of Fibonacci elements in Coding Theory. The matrix given in
relation (3.9) can be used for coding and decoding over $\mathbb{Z}$ when
its determinant is $1\,$or $-1$ or for coding and decoding over $\mathbb{Z}%
_{a}$ when its determinant is prime with the integer $a$. Even if these
matrices have the same properties for coding and decoding as the similar
matrices generated by Fibonacci elements, that means can correct $1,2$ or $3$
errors, they can be used as a new class of such a matrices utilized for
coding and decoding.\medskip

\begin{equation*}
\end{equation*}%
\textbf{4. Generalized Fibonacci- Lucas quaternions. Some properties and
applications}

\begin{equation*}
\end{equation*}

\bigskip

Let $\gamma _{1},\gamma _{2}\in \mathbb{R}\setminus \{0\}$ and let $\mathbb{H%
}\left( \gamma _{1},\gamma _{2}\right) $ be the generalized quaternion%
{\small \ }algebra with basis $\{1,e_{1},e_{2},e_{3}\},$ that means a real
algebra with the multiplication given in the following table%
\begin{equation*}
\begin{tabular}{c||c|c|c|c|}
$\cdot $ & $1$ & $e_{1}$ & $e_{2}$ & $e_{3}$ \\ \hline\hline
$1$ & $1$ & $e_{1}$ & $e_{2}$ & $e_{3}$ \\ \hline
$e_{1}$ & $e_{1}$ & $\gamma _{1}$ & $e_{3}$ & $\gamma _{1}e_{2}$ \\ \hline
$e_{2}$ & $e_{2}$ & $-e_{3}$ & $\gamma _{2}$ & $-\gamma _{2}e_{1}$ \\ \hline
$e_{3}$ & $e_{3}$ & $-\gamma _{1}e_{2}$ & $\gamma _{2}e_{1}$ & $\gamma
_{1}\gamma _{2}$ \\ \hline
\end{tabular}%
.
\end{equation*}%
We consider $\alpha ,\beta $$\in $$\mathbb{Q}^{\ast }$ and $\mathbb{H}_{%
\mathbb{Q}}\left( \alpha ,\beta \right) $ the generalized quaternion algebra
over the rational field. Let $d_{n}$ be $(a,b,x_{0},x_{1})-$numbers. We \
define the $n$-\textit{th}\newline
$(a,b,x_{0},x_{1})-$\textit{quaternions} to be the elements of the form 
\begin{equation*}
D_{n}=d_{n}+d_{n+1}e_{1}+d_{n+2}e_{2}+d_{n+3}e_{3}
\end{equation*}

In the paper [Fl, Sa; 15], we introduced the $n$-\textit{th generalized
Fibonacci-Lucas quaternion} to be the elements of the form 
\begin{equation*}
G_{n}^{p,q}=g_{n}^{p,q}1+g_{n+1}^{p,q}e_{1}+g_{n+2}^{p,q}e_{2}+g_{n+3}^{p,q}e_{3}.
\end{equation*}%
We remark that these elements are $n$-\textit{th} $\left( 1,1,p+2q,q\right)
- $quaternions.\medskip

\textbf{Remark 4.1.} In the paper [Fl, Sa; 15], Theorem 3.5, we proved that
the set 
\begin{equation*}
\left\{ \sum\limits_{i=1}^{n}5G_{n_{i}}^{p_{i},q_{i}}|n\in \mathbb{N}^{\ast
},p_{i},q_{i}\in \mathbb{Z},(\forall )i=\overline{1,n}\right\} \cup \left\{
1\right\}
\end{equation*}%
has a ring structure with quaternions addition and multiplication, that
means it is an order of a rational quaternion algebra. To prove that the
above set is closed under multiplications, we used properties of Fibonacci
and Lucas number regarding multiplications of two such elements (see
Proposition 2.1). For example, we remark that the product of two Fibonacci
numbers gives us Lucas numbers. This remark explains why we considered the
generalized Fibonacci-Lucas numbers given by the relation 1.2.

In the following, we will provide some properties of generalized
Fibonacci-Lucas quaternions.

Let $B$ be the generating function for the generalized Fibonacci-Lucas
quaternions, $B\left( z\right) =\sum\limits_{n\geq
1}G_{n}^{p,q}z^{n}.\medskip $

\textbf{Proposition 4.2.} \textit{With the above notations, we have} 
\begin{equation*}
B\left( z\right) =\frac{G_{1}^{p,q}z+\left( G_{2}^{p,q}-G_{1}^{p,q}\right)
z^{2}}{1-z-z^{2}}.
\end{equation*}%
\textbf{Proof.} We use the following relations\newline
\begin{equation}
B\left( z\right) =G_{1}^{p,q}z+G_{2}^{p,q}z^{2}+...G_{n}^{p,q}z^{n}+...\ \  
\tag{4.1}
\end{equation}%
\begin{equation}
zB\left( z\right)
=G_{1}^{p,q}z^{2}+G_{2}^{p,q}z^{3}+...G_{n-1}^{p,q}z^{n}+...\ \   \tag{4.2}
\end{equation}%
\begin{equation}
z^{2}B\left( z\right)
=G_{1}^{p,q}z^{3}+G_{2}^{p,q}z^{4}+...G_{n-2}^{p,q}z^{n}+...\ \ .  \tag{4.3}
\end{equation}%
Adding the equalities (4.1) and (4.2) member by member and using relation
(4.3), we obtain 
\begin{equation*}
B\left( z\right) \left( 1-z-z^{2}\right) =G_{1}^{p,q}z+\left(
G_{2}^{p,q}-g_{1}^{p,q}\right) z^{2},
\end{equation*}%
therefore%
\begin{equation*}
B\left( z\right) =\frac{G_{1}^{p,q}\cdot z+\left(
G_{2}^{p,q}-g_{1}^{p,q}\right) \cdot z^{2}}{1-z-z^{2}}.
\end{equation*}%
$\Box \smallskip $\medskip

In the following, we consider $h(x)$ a polynomial with real coefficients. We
define the real $h(x)$--\textit{generalized Fibonacci-Lucas polynomials}, as
polynomials defined by the recurrence relation 
\begin{equation}
g_{h,n}^{p,q}(x)=h(x)g_{h,n-1}^{p,q}(x)+g_{h,n-2}^{p,q}(x),\quad
n=2,3,\ldots ,  \tag{4.4}
\end{equation}%
where $g_{h,0}^{p,q}(x)=p+2q,$ $g_{h,1}^{p,q}(x)=q,$ for all polynomials $%
h(x)$.\medskip

\textbf{Definition 4.3. }The $h(x)$ -- \textit{generalized quaternion
Fibonacci-Lucas polynomials} $\{G_{h,n}^{p,q}(x)\}_{n=0}^{\infty }$ are
given by the recurrence relation 
\begin{equation}
G_{h,n}^{p,q}(x)=\sum\limits_{k=0}^{3}g_{h,n+k}^{p,q}(x)\,e_{k},  \tag{4.5}
\end{equation}%
where $g_{h,n}^{p,q}(x)$ is the $n$th real $h(x)$ -- Fibonacci-Lucas
polynomial.\medskip

\textbf{Definition 4.4.} The \textit{generating function} $A(t)$
corresponding to the sequence $\{G_{h,n}^{p,q}(x)\}_{n=0}^{\infty }$ is
defined by the following relation 
\begin{equation}
A(t)=\sum\limits_{n=0}^{\infty }G_{h,n}^{p,q}(x)t^{n}.  \tag{4.6}
\end{equation}

\textbf{Theorem 4.5.} \textit{The generating function for } $h(x)$ \textit{%
-- quaternion Fibonacci-Lucas polynomials} $G_{h,n}^{p,q}(x)$ \textit{is
given by the relation} 
\begin{equation*}
A(t)=\frac{G_{h,0}^{p,q}(x)+(G_{h,1}^{p,q}(x)-h(x)G_{h,0}^{p,q}(x))t}{%
1-h(x)t-t^{2}}.
\end{equation*}

\textbf{Proof. }From relation (4.6), we obtain the following relation 
\begin{equation*}
A(t)(1-h(x)t-t^{2})=
\end{equation*}%
\begin{equation*}
=\sum\limits_{n=0}^{\infty
}G_{h,n}^{p,q}(x)t^{n}-h(x)\sum\limits_{n=0}^{\infty
}G_{h,n}^{p,q}(x)t^{n+1}-\sum\limits_{n=0}^{\infty }G_{h,n}^{p,q}(x)t^{n+2}=
\end{equation*}%
\begin{equation*}
=G_{h,0}^{p,q}(x)+(G_{h,1}^{p,q}(x)-h(x)G_{h,0}^{p,q}(x))t+
\end{equation*}%
\begin{equation*}
+\sum\limits_{n=2}^{\infty
}t^{n}(G_{h,n}^{p,q}-h(x)G_{h,n-1}^{p,q}-G_{h,n-2}^{p,q})=
\end{equation*}%
\begin{equation*}
=G_{h,0}(x)+(G_{h,1}(x)-h(x)G_{h,0}(x))t.
\end{equation*}%
$\Box \smallskip $

Let $r_{1}(x)$ and $r_{2}(x)$ be the solutions of the characteristic
equation $r^{2}-h(x)r-1=0,$ for the recurrence relation given in \ (4.4). We
obtain that 
\begin{equation}
r_{1}(x)=\frac{h(x)+\sqrt{h^{2}(x)+4}}{2}\,,\quad r_{2}(x)=\frac{h(x)-\sqrt{%
h^{2}(x)+4}}{2}\,.  \tag{4.7}
\end{equation}%
\qquad \qquad

\textbf{Theorem 4.6.} (Binet's formula) \textit{For all} $n\in
\{0,1,2,\ldots \},$ \textit{we have the following relation for the
polynomials} $g_{h,n}^{p,q}(x)$ 
\begin{equation}
g_{h,n}^{p,q}(x)=\frac{\left( p+2q\right) }{r_{1}(x)-r_{2}(x)}\left(
r_{1}^{n+1}(x)-r_{2}^{n+1}(x)\right) -\frac{q}{r_{1}(x)-r_{2}(x)}\left(
r_{1}^{n}(x)-r_{2}^{n}(x)\right)  \tag{4.8}
\end{equation}

\textbf{Proof.} Assuming that $g_{h,n}^{p,q}(x)=\alpha r_{1}^{n}(x)+\beta
r_{2}^{n}(x).$ From relation (4.4), we obtain that $\alpha +\beta =$ $p+2q$
and $\alpha r_{1}(x)+\beta r_{2}(x)=$ $q.$ Solving this system, we obtain 
\begin{equation*}
\alpha =\frac{\left( p+2q\right) r_{1}\left( x\right) -q}{r_{1}(x)-r_{2}(x)}
\end{equation*}%
and 
\begin{equation*}
\beta =\frac{q-\left( p+2q\right) r_{2}\left( x\right) }{r_{1}(x)-r_{2}(x)}.
\end{equation*}%
Therefore 
\begin{equation*}
g_{h,n}^{p,q}(x)=\frac{\left( p+2q\right) r_{1}\left( x\right) -q}{%
r_{1}(x)-r_{2}(x)}r_{1}^{n}(x)+\frac{q-\left( p+2q\right) r_{2}\left(
x\right) }{r_{1}(x)-r_{2}(x)}r_{2}^{n}(x).
\end{equation*}%
It results \ that $g_{h,n}^{p,q}(x)=\frac{\left( p+2q\right) }{%
r_{1}(x)-r_{2}(x)}\left( r_{1}^{n+1}(x)-r_{2}^{n+1}(x)\right) -\frac{q}{%
r_{1}(x)-r_{2}(x)}\left( r_{1}^{n}(x)-r_{2}^{n}(x)\right) .\Box \smallskip $

\textbf{Theorem 4.7.} \textit{For} $n\in \{0,1,2,\ldots \},$ \textit{we have
the following relation} 
\begin{equation}
G_{h,n}^{p,q}(x)=\frac{R_{1}(x)r_{1}^{n}(x)-R_{2}(x)r_{2}^{n}(x)}{%
r_{1}(x)-r_{2}(x)},  \tag{4.9}
\end{equation}%
\textit{where}\newline
$R_{1}\left( x\right) =\sum\limits_{k=0}^{3}\left( -qr_{1}^{k}\left(
x\right) +\left( p+2q\right) r_{1}^{k+1}\left( x\right) \right)
\,e_{k}\,,\quad $\newline
$R_{2}\left( x\right) =\sum\limits_{k=0}^{3}\left( -qr_{2}^{k}\left(
x\right) +\left( p+2q\right) r_{2}^{k+1}\left( x\right) \right) \,e_{k}\,\,.$

\textbf{Proof.} Using the Binet's formula (4.8), it results\newline
$G_{h,n}^{p,q}(x)=\sum\limits_{k=0}^{3}g_{h,n+k}^{p,q}(x)\,e_{k}=$\newline
$=\frac{\left( p+2q\right) }{r_{1}(x)-r_{2}(x)}\left(
r_{1}^{n+1}(x)-r_{2}^{n+1}(x)\right) e_{0}-\frac{q}{r_{1}(x)-r_{2}(x)}\left(
r_{1}^{n}(x)-r_{2}^{n}(x)\right) e_{0}+$\newline
$+\frac{\left( p+2q\right) }{r_{1}(x)-r_{2}(x)}\left(
r_{1}^{n+2}(x)-r_{2}^{n+2}(x)\right) e_{1}-\frac{q}{r_{1}(x)-r_{2}(x)}\left(
r_{1}^{n+1}(x)-r_{2}^{n+1}(x)\right) e_{1}+$\newline
$+\frac{\left( p+2q\right) }{r_{1}(x)-r_{2}(x)}\left(
r_{1}^{n+3}(x)-r_{2}^{n+3}(x)\right) e_{2}-\frac{q}{r_{1}(x)-r_{2}(x)}\left(
r_{1}^{n+2}(x)-r_{2}^{n+2}(x)\right) e_{2}+$\newline
$+\frac{\left( p+2q\right) }{r_{1}(x)-r_{2}(x)}\left(
r_{1}^{n+4}(x)-r_{2}^{n+4}(x)\right) e_{3}-\frac{q}{r_{1}(x)-r_{2}(x)}\left(
r_{1}^{n+3}(x)-r_{2}^{n+3}(x)\right) e_{3}=$\newline
$=\frac{r_{1}^{n}(x)}{r_{1}(x)-r_{2}(x)}[-q+\left( p+2q\right) r_{1}\left(
x\right) +\left( -qr_{1}\left( x\right) +\left( p+2q\right) r_{1}^{2}\left(
x\right) \right) e_{1}+$\newline
$+\left( -qr_{1}^{2}\left( x\right) +\left( p+2q\right) r_{1}^{3}\left(
x\right) \right) e_{2}+\left( -qr_{1}^{3}\left( x\right) +\left( p+2q\right)
r_{1}^{4}\left( x\right) \right) e_{3}]-$\newline
$-\frac{r_{2}^{n}(x)}{r_{1}(x)-r_{2}(x)}[-q+\left( p+2q\right) r_{2}\left(
x\right) +\left( -qr_{2}\left( x\right) +\left( p+2q\right) r_{2}^{2}\left(
x\right) \right) e_{1}+$\newline
$+\left( -qr_{2}^{2}\left( x\right) +\left( p+2q\right) r_{2}^{3}\left(
x\right) \right) e_{2}+\left( -qr_{2}^{3}\left( x\right) +\left( p+2q\right)
r_{2}^{4}\left( x\right) \right) e_{3}]=$\newline
$=R_{1}\left( x\right) \frac{r_{1}^{n}(x)}{r_{1}(x)-r_{2}(x)}+R_{2}\left(
x\right) \frac{r_{2}^{n}(x)}{r_{1}(x)-r_{2}(x)}$%
\begin{equation*}
=\frac{R_{1}(x)r_{1}^{n}(x)-R_{2}(x)r_{2}^{n}(x)}{r_{1}(x)-r_{2}(x)}\,.
\end{equation*}%
$\Box \smallskip $\medskip

\textbf{Theorem 4.8.} (Catalan's identity) \textit{Let }$n,s$ \textit{be
positive integers} \textit{with} $s\leq n$. \textit{The following relation
is true} 
\begin{equation*}
G_{h,n+s}^{p,q}(x)G_{h,n-s}^{p,q}(x)-G_{h,n}^{p,q2}(x)=
\end{equation*}%
\begin{equation*}
\frac{(-1)^{n+s+1}}{h^{2}(x)+4}%
[R_{1}(x)R_{2}(x)((-1)^{s+1}+r_{1}^{2}(x))+R_{1}(x)R_{2}(x)((-1)^{s+1}+r_{2}^{2}(x))].
\end{equation*}

\textbf{Proof.} Using formula (4.9), it results 
\begin{equation*}
G_{h,n+s}^{p,q}(x)G_{h,n-s}^{p,q}(x)-G_{h,n}^{p,q2}(x)=
\end{equation*}%
\begin{equation*}
\frac{1}{(r_{1}(x)-r_{2}(x))^{2}}%
[R_{1}(x)R_{2}(x)r_{1}^{n}(x)r_{2}^{n}(x)(1-\left( \frac{r_{1}(x)}{r_{2}(x)}%
\right) ^{s}\,)+
\end{equation*}%
\begin{equation*}
R_{1}(x)R_{2}(x)r_{1}^{n}(x)r_{2}^{n}(x)(1-\left( \frac{r_{2}\left( x\right) 
}{r_{1}(x)}\right) ^{r}\,)].
\end{equation*}%
Finally, we use the Vi\`{e}te's relations, we obtain the asked identity.$%
\Box \smallskip $

If in the above Theorem , we take $r=1$, we obtain the Cassini's
identity.\medskip

\textbf{Theorem 4.9} (Cassini's identity) \textit{For each natural number }$%
n $\textit{, we have} 
\begin{equation*}
G_{h,n+1}^{p,q}(x)G_{h,n-1}^{p,q}(x)-G_{h,n}^{p,q2}(x)=
\end{equation*}%
\begin{equation*}
\frac{(-1)^{n}}{h^{2}(x)+4}%
[R_{1}(x)R_{2}(x)(1+r_{1}^{2}(x))+R_{1}(x)R_{2}(x)(1+r_{2}^{2}(x))].
\end{equation*}%
$\Box \smallskip $

Similar results with the results obtained above for the generalized
Fibonacci-Lucas quaternions were obtained in \textbf{[}Ca; 15\textbf{], }%
Theorem 3.3 and Theorem 3.6 for $h\left( x\right) -$Fibonacci polynomials
over the real field and in \textbf{[}Fl, Sh, Vl; 17\textbf{], }Theorem 2.3
and Theorem 2.6, for $h\left( x\right) -$Fibonacci polynomials over an
arbitrary algebra.\medskip 
\begin{equation*}
\end{equation*}%
\textbf{5. Some properties of } $\left( 1,a,0,1\right) -$\textbf{quaternions
and} $\left( 1,a,2,1\right) -$\textbf{quaternions}%
\begin{equation*}
\end{equation*}

Let $a$ be a nonzero natural number. Let $\left( x_{n}\right) _{n\geq 0}$ be
the $\left( 1,a,0,1\right) -$numbers and $\left( y_{n}\right) _{n\geq 0}$ be
the $\left( 1,a,2,1\right) -$numbers, that means 
\begin{equation*}
x_{n}=x_{n-1}+ax_{n-2},\;n\geq 2,x_{0}=0,x_{1}=1
\end{equation*}%
and 
\begin{equation*}
y_{n}=y_{n-1}+ay_{n-2},\;n\geq 2,y_{0}=2,y_{1}=1.
\end{equation*}%
If $n$ is a negative integer, we take $x_{n}=\left(-1\right)^{-n+1} \cdot
x_{-n}.$\newline
In the paper [Fl, Sa; 15], we introduced the generalized Fibonacci-Lucas
numbers and the generalized Fibonacci-Lucas quaternions and we obtained some
properties of them.\newline
Now, we introduce another numbers and another quaternions, where instead of
the Fibonacci sequence $\left( f_{n}\right) _{n\geq 0}$ we consider the
sequence $\left( x_{n}\right) _{n\geq 0}$ and instead of the Lucas sequence $%
\left( l_{n}\right) _{n\geq 0}$ we consider the sequence $\left(
y_{n}\right) _{n\geq 0}.$\newline
If we denote with $\alpha =\frac{1+\sqrt{1+4a}}{2}$ and $\beta =\frac{1-%
\sqrt{1-4a}}{2},$ it is easy to obtain the following relations: \smallskip 
\newline
\textbf{Binet's formula for the sequence} $\left( x_{n}\right) _{n\geq 0}.$ 
\begin{equation*}
x_{n}=\frac{\alpha ^{n}-\beta ^{n}}{\alpha -\beta }=\frac{\alpha ^{n}-\beta
^{n}}{\sqrt{1+4a}},\ \ \left( \forall \right) n\in \mathbb{N}.
\end{equation*}%
\textbf{Binet's formula for the sequence} $\left( y_{n}\right) _{n\geq 0}.$ 
\begin{equation*}
y_{n}=\alpha ^{n}+\beta ^{n},\ \ \left( \forall \right) n\in \mathbb{N}.
\end{equation*}%
First of all, we get some properties of the sequences $\left( x_{n}\right)
_{n\geq 0},$ $\left( y_{n}\right) _{n\geq 0}.$\newline
\smallskip\newline
\textbf{Proposition 5.1.} \textit{Let} $\left( x_{n}\right) _{n\geq 0},$ $%
\left( y_{n}\right) _{n\geq 0}$ \textit{be the sequences previously defined.
Then, we have}:\newline
i) 
\begin{equation*}
y_{n}y_{n+l}=y_{2n+l}+\left( -a\right) ^{n}y_{l},\ \left( \forall \right)
n,l\in \mathbb{N};
\end{equation*}%
ii) 
\begin{equation*}
x_{n}y_{n+l}=x_{2n+l}-\left( -a\right) ^{n}x_{l},\ \left( \forall \right)
n,l\in \mathbb{N};
\end{equation*}%
iii) 
\begin{equation*}
x_{n+l}y_{n}=x_{2n+l}+\left( -a\right) ^{n}x_{l},\ \left( \forall \right)
n,l\in \mathbb{N};
\end{equation*}%
iv) 
\begin{equation*}
x_{n}x_{n+l}=\frac{1}{1+4a}\left[ y_{2n+l}-\left( -a\right) ^{n}y_{l}\right]
,\ \left( \forall \right) n,l\in \mathbb{N}.
\end{equation*}%
\textbf{Proof.} Using Binet's formulae for the sequences $\left(
x_{n}\right) _{n\geq 0}$ and $\left( y_{n}\right) _{n\geq 0},$ we have:%
\newline
i) 
\begin{equation*}
y_{n}y_{n+l}=\left( \alpha ^{n}+\beta ^{n}\right) \left( \alpha ^{n+l}+\beta
^{n+l}\right) =\alpha ^{2n+l}+\alpha ^{n}\beta ^{n+l}+\alpha ^{n+l}\beta
^{n}+\beta ^{2n+l}=
\end{equation*}%
\begin{equation*}
=y_{2n+l}+\alpha ^{n}\beta ^{n}\left( \alpha ^{l}+\beta ^{l}\right)
=y_{2n+l}+\left( -a\right) ^{n}y_{l}.
\end{equation*}%
ii) 
\begin{equation*}
x_{n}y_{n+l}=\frac{\alpha ^{n}-\beta ^{n}}{\alpha -\beta }\left( \alpha
^{n+l}+\beta ^{n+l}\right) =
\end{equation*}%
\begin{equation*}
=\frac{\alpha ^{2n+l}-\beta ^{2n+l}}{\alpha -\beta }-\frac{\alpha ^{n}\beta
^{n}\left( \alpha ^{l}-\beta ^{l}\right) }{\alpha -\beta }=x_{2n+l}-\left(
-a\right) ^{n}x_{l}.
\end{equation*}%
iii) 
\begin{equation*}
x_{n+l}y_{n}=\frac{\alpha ^{n+l}-\beta ^{n+l}}{\alpha -\beta }\left( \alpha
^{n}+\beta ^{n}\right) =
\end{equation*}%
\begin{equation*}
=\frac{\alpha ^{2n+l}-\beta ^{2n+l}}{\alpha -\beta }+\frac{\alpha ^{n}\beta
^{n}\left( \alpha ^{l}-\beta ^{l}\right) }{\alpha -\beta }=x_{2n+l}+\left(
-a\right) ^{n}x_{l}.
\end{equation*}%
iv) 
\begin{equation*}
x_{n}x_{n+l}=\frac{\alpha ^{n}-\beta ^{n}}{\alpha -\beta }\cdot \frac{\alpha
^{n+l}-\beta ^{n+l}}{\alpha -\beta }=
\end{equation*}%
\begin{equation*}
=\frac{\alpha ^{2n+l}+\beta ^{2n+l}-\alpha ^{n}\beta ^{n}\left( \alpha
^{l}+\beta ^{l}\right) }{\left( \alpha -\beta \right) ^{2}}=\frac{1}{1+4a}%
\left[ y_{2n+l}-\left( -a\right) ^{n}y_{l}\right] .
\end{equation*}%
$\Box \smallskip $\medskip \newline
Let $p,q$ be two arbitrary integers. We consider the sequence $\left(
s_{n}\right) _{n\geq 0},$ 
\begin{equation}
s_{n+1}=px_{n}+qy_{n+1},\;n\geq 0,  \tag{5.1.}
\end{equation}%
where $\left( x_{n}\right) _{n\geq 0},$ $\left( y_{n}\right) _{n\geq 0}$ are
the sequences previously defined.\newline
We obtain that $s_{n}=s_{n-1}+as_{n-2},\ \left( \forall \right) n\in \mathbb{%
N},n\geq 2,$\newline
$s_{0}=px_{-1}+qy_{0}=p+2q,$ $s_{1}=q,$ that means $\left( s_{n}\right)
_{n\geq 0}$ are $\left( 1,a,p+2q,q\right) -$numbers.\newline
In the following, we will use the notation $s_{n}^{p,q}$ for $s_{n}.$\newline

Let $\alpha ,\beta $$\in $$\mathbb{Q}^{\ast }$ and let $\mathbb{H}_{\mathbb{Q%
}}\left( \alpha ,\beta \right) $ be the generalized quaternion algebra over
the rational field, with basis $\{1,e_{1},e_{2},e_{3}\}.$ We define the $n$%
-th\newline
$\left( 1,a,p+2q,q\right)-$ quaternions to be the elements of the form 
\begin{equation*}
S_{n}^{p,q}=s_{n}^{p,q}1+s_{n+1}^{p,q}e_{1}+s_{n+2}^{p,q}e_{2}+s_{n+3}^{p,q}e_{3}.
\end{equation*}%
\textbf{Remark 5.2.} \textit{Let} $p,q$ \textit{be two arbitrary integers,
let} $n$ \textit{be an arbitrary positive integer and let} $\left(
s_{n}^{p,q}\right) _{n\geq 1}$ \textit{the sequence previously defined.
Then, we have:} 
\begin{equation*}
px_{n+1}+qy_{n}=s_{n}^{ap,q}+s_{n+1}^{p,0},\forall ~n\in \mathbb{N}-\{0\}.
\end{equation*}%
\textbf{Proof.} We compute 
\begin{equation*}
px_{n+1}+qy_{n}=px_{n}+apx_{n-1}+qy_{n}=s_{n}^{ap,q}+s_{n+1}^{p,o}.
\end{equation*}%
$\Box \smallskip $\medskip \newline
\textbf{Remark 5.3. } \textit{Using the previously notations, we have the
following relation} 
\begin{equation*}
S_{n}^{p,q}=0\ \text{\textit{if\ and\ only\ if}}\ p=q=0.
\end{equation*}%
\textbf{Proof.} " $\Rightarrow $" If $S_{n}^{p,q}=0,$ since $%
\{1,e_{1},e_{2},e_{3}\}$ is a $\mathbb{Q}-$ basis in quaternion algebra $%
\mathbb{H}_{\mathbb{Q}}\left( \alpha ,\beta \right) ,$ it results that $%
s_{n}^{p,q}=s_{n+1}^{p,q}=s_{n+2}^{p,q}=s_{n+3}^{p,q}=0.$ Using the
recurrence relation of the sequence $\left( s_{n}^{p,q}\right) _{n\geq 1},$
we obtain that $s_{n-1}^{p,q}=0,$ $s_{n-2}^{p,q}=0,$ ..., $s_{1}^{p,q}=0$, $%
s_{0}^{p,q}=0.$ Since $s_{1}^{p,q}=q,$ it results $q=0$ and since $%
s_{0}^{p,q}=p+2q=0,$ it results $p=0.$\newline
" $\Leftarrow $" It is trivial. $\Box \smallskip $\medskip \newline
\textbf{Proposition 5.4.} \textit{Let} $a$ \textit{be a nonzero natural
number and let} $O$ \textit{be} \textit{the set} 
\begin{equation*}
O=\left\{ \sum\limits_{i=1}^{n}\left( 1+4a\right)
S_{n_{i}}^{p_{i},q_{i}}|n\in \mathbb{N}^{\ast },p_{i},q_{i}\in \mathbb{Z}%
,(\forall )i=\overline{1,n}\right\} \cup \left\{ 1\right\} .
\end{equation*}%
\textit{Then} $O$ \textit{is an order of the quaternion algebra} $\mathbb{H}%
_{\mathbb{Q}}\left( \alpha ,\beta \right) .$ \newline
\smallskip \newline
\textbf{Proof.} Applying Remark 5.3, $S_{n}^{0,0}=0$$\in $$O.$\newline
Let $n,m\in \mathbb{N}^{\ast },$ $p,q,p^{^{\prime }},q^{^{\prime }},c,d\in 
\mathbb{Z}.$ We obtain that 
\begin{equation*}
cs_{n}^{p,q}+ds_{m}^{p^{^{\prime }},q^{^{\prime
}}}=s_{n}^{cp,cq}+s_{m}^{dp^{^{\prime }},dq^{^{\prime }}}
\end{equation*}%
and, from here, we get 
\begin{equation*}
cS_{n}^{p,q}+dS_{m}^{p^{^{\prime }},q^{^{\prime
}}}=S_{n}^{cp,cq}+S_{m}^{dp^{^{\prime }},dq^{^{\prime }}}.
\end{equation*}%
This implies that $O$ \textit{is a free} $\mathbb{Z}-$ submodule of rank $4$
of the quaternion algebra $\mathbb{H}_{\mathbb{Q}}\left( \alpha ,\beta
\right) .$ \newline
Now, we prove that $O$ is a subring of $\mathbb{H}_{\mathbb{Q}}\left( \alpha
,\beta \right) .$ Let $m,n$ be two integers, $n<m.$ We have:\newline
\begin{equation*}
\left( 1+4a\right) s_{n}^{p,q}\left( 1+4a\right) s_{m}^{p^{^{\prime
}},q^{^{\prime }}}\text{=}\left( 1+4a\right) \left( px_{n-1}+qy_{n}\right)
\left( 1+4a\right) \left( p^{^{\prime }}x_{m-1}+q^{^{\prime }}y_{m}\right) 
\text{=}
\end{equation*}%
\begin{equation*}
=\left( 1+4a\right) ^{2}pp^{^{\prime }}x_{n-1}x_{m-1}+\left( 1+4a\right)
^{2}pq^{^{\prime }}x_{n-1}y_{m}+
\end{equation*}%
\begin{equation*}
+\left( 1+4a\right) ^{2}p^{^{\prime }}qx_{m-1}y_{n}+\left( 1+4a\right)
^{2}qq^{^{\prime }}y_{n}y_{m}.
\end{equation*}%
Applying Proposition 5.1 and Remark 5.2, we obtain: 
\begin{equation*}
\left( 1+4a\right) s_{n}^{p,q}\left( 1+4a\right) s_{m}^{p^{^{\prime
}},q^{^{\prime }}}\text{=}\left( 1+4a\right) ^{2}pp^{^{\prime }}\frac{1}{1+4a%
}\left[ y_{n+m-2}-\left( -a\right) ^{n}y_{m-n}\right] +
\end{equation*}%
\begin{equation*}
+\left( 1+4a\right) ^{2}pq^{^{\prime }}\left[ x_{m+n-1}-\left( -a\right)
^{n-1}x_{m-n+1}\right] +
\end{equation*}%
\begin{equation*}
\text{+}\left( 1+4a\right) ^{2}p^{^{\prime }}q\left[ x_{n+m-1}+\left(
-a\right) ^{n}x_{m-n-1}\right] \text{+}\left( 1+4a\right) ^{2}qq^{^{\prime }}%
\left[ y_{n+m}+\left( -a\right) ^{n}y_{m-n}\right] \text{=}
\end{equation*}%
\begin{equation*}
\text{=}\left( 1\text{+}4a\right) ^{2}\left[ pq^{^{\prime }}x_{n+m-1}\text{+}%
qq^{^{\prime }}y_{m+n}\right] \text{+}\left( 1\text{+}4a\right) ^{2}\left[
\left( \text{-}a\right) ^{n}p^{^{\prime }}qx_{m-n-1}\text{+}\left( \text{-}%
a\right) ^{n}qq^{^{\prime }}y_{m-n}\right] \text{+}
\end{equation*}%
\begin{equation*}
+\left( 1+4a\right) \left[ -\left( -a\right) ^{n-1}\left( 1+4a\right)
pq^{^{\prime }}x_{m-n+1}-\left( -a\right) ^{n}pp^{^{\prime }}y_{m-n}\right] +
\end{equation*}%
\begin{equation*}
+\left( 1+4a\right) \left[ \left( 1+4a\right) p^{^{\prime
}}qx_{n+m-1}+pp^{^{\prime }}y_{n+m-2}\right] =
\end{equation*}%
\begin{equation*}
\text{=}\left( 1+4a\right) s_{m+n}^{\left( 1+4a\right) pq^{^{\prime
}},\left( 1+4a\right) qq^{^{\prime }}}\text{+}\left( 1+4a\right)
s_{m-n}^{\left( -a\right) ^{n}\left( 1+4a\right) p^{^{\prime }}q,\left(
-a\right) ^{n}\left( 1+4a\right) qq^{^{\prime }}}\text{+}
\end{equation*}%
\begin{equation*}
+\left( 1+4a\right) s_{m-n+1}^{\left( -a\right) ^{n}\left( 1+4a\right)
pq^{^{\prime }},\left( -a\right) ^{n+1}pp^{^{\prime }}}+\left( 1+4a\right)
s_{m-n+1}^{\left( -a\right) ^{n}\left( 1+4a\right) pq^{^{\prime }},0}+
\end{equation*}%
\begin{equation*}
+\left( 1+4a\right) s_{m+n-2}^{a\left( 1+4a\right) p^{^{\prime
}}q,pp^{^{\prime }}}+\left( 1+4a\right) s_{m+n-1}^{a\left( 1+4a\right)
p^{^{\prime }}q,0}.
\end{equation*}%
It results that $\left( 1+4a\right) s_{n}^{p,q}\left( 1+4a\right)
s_{m}^{p^{^{\prime }},q^{^{\prime }}}$$\in $$O.$ Therefore, $O$ is an order
of the quaternion algebra $\mathbb{H}_{\mathbb{Q}}\left( \alpha ,\beta
\right) .\Box \smallskip $\medskip \newline
\textbf{Conclusions.} In this paper, we introduced the $(a,b,x_{0},x_{1})-$%
elements and we studied some properties and applications for the $\left(
1,1,p+2q,q\right) -$numbers (the generalized Fibonacci-Lucas numbers), the $%
\left( 1,1,p+2q,q\right) -$quaternions (generalized Fibonacci-Lucas
quaternions), the $\left( 1,a,0,1\right) -$numbers, the \newline
$\left( 1,a,2,1\right) -$numbers and the $\left( 1,a,p+2q,q\right)-$%
quaternions. For the last one, we gave an interesting algebraic structure.

From the above, we can see that these elements can be studied \ in a more
genealized cases and all other particular cases are subordinated to this
approach. In the paper [Fl, Sa; 17], we studied properties and applications 
of elements arising from a difference equation of degree three. It is interesting 
to study what relations must satisfy  coefficients of a difference equation such that 
the quaternions defined using that equation to have  a ring structure, as above. This idea
can constitute a starting point for a further research.

 \newpage 
\begin{equation*}
\end{equation*}%
\textbf{References}\newline
\begin{equation*}
\end{equation*}%
\newline
[Ak, Ko, To; 14]\textbf{\ }M. Akyigit, H. H. Koksal, M. Tosun, \textit{%
Fibonacci Generalized Quaternions}, Adv. Appl. Clifford Algebras,
3(24)(2014), 631--641.\newline
[Ba, Pr; 09] M. Basu, B. Prasad, \textit{The generalized relations among the
code elements for Fibonacci coding theory}, Chaos, Solitons and Fractals,
41(2009), 2517-2525.\newline
[Ca; 15] P. Catarino, \textit{A note on} $h(x)$ \textit{-- fibonacci
quaternion polynomials}, Chaos, Solitons and Fractals, 77(2015), 1--5.%
\newline
[Fa, Pl; 07(1)] S. Falc\'{o}n, \'{A}. Plaza, \textit{On the Fibonacci} $k$%
\textit{-numbers}, Chaos, Solitons and Fractals, 32(5)(2007), 1615--24.%
\newline
[Fa, Pl; 07(2)] S. Falc\'{o}n, \'{A}. Plaza, \textit{The} $k$\textit{%
-Fibonacci sequence and the Pascal 2-triangle}, Chaos, Solitons and
Fractals, 33(1)(2007), 38--49.\newline
[Fa, Pl; 09] S. Falc\'{o}n, A. Plaza, \textit{On} $k$\textit{-Fibonacci
sequences and polynomials and their derivatives}, Chaos, Solitons and
Fractals 39(3)(2009), 1005--1019.\newline
[Fib.] http://www.maths.surrey.ac.uk/hosted-sites/R.Knott/Fibonacci/fib.html%
\newline
[Fl, Sa; 15] C. Flaut, D. Savin, \textit{Quaternion Algebras and Generalized
Fibonacci-Lucas Quaternions}, Adv. Appl. Clifford Algebras, 25(4)(2015),
853-862.\newline
[Fl, Sa; 17] C. Flaut, D. Savin, \textit{Applications of some special numbers obtained from a
difference equation of degree three}, arXiv:1705.02677.\newline
[Fl, Sh;13] C. Flaut, V. Shpakivskyi, \textit{On Generalized Fibonacci
Quaternions and Fibonacci-Narayana Quaternions}, Adv. Appl. Clifford
Algebras, 23(3)(2013), 673--688.\newline
[Fl, Sh; 15] C. Flaut, V. Shpakivskyi, \textit{Some remarks about Fibonacci
elements in an arbitrary algebra}, Bull. Soc. Sci. Lettres \L \'{o}d\'{z},
65(3)(2015), 63--73.\newline
[Fl, Sh, Vl; 17] C. Flaut, V. Shpakivskyi, E. Vlad, \textit{Some remarks
regarding} $h(x)$ \textit{-- Fibonacci polynomials in an arbitrary algebra},
Chaos, Solitons \& Fractals, 99(2017), 32-35.\newline
[Gu, Nu; 15]\ I. A. Guren, S.K. Nurkan, \textit{A new approach to Fibonacci,
Lucas numbers and dual vectors}, Adv. Appl. Clifford Algebras, 3(25)(2015),
577--590.\newline
[Ha; 12]\ S. Halici, \textit{On Fibonacci Quaternions}, Adv. in Appl.
Clifford Algebras, 22(2)(2012), 321--327.\newline
[Ho; 63]\ A. F. Horadam, \textit{Complex Fibonacci Numbers and Fibonacci
Quaternions}, Amer. Math. Monthly, 70(1963), 289--291.\newline
[Na, Ha; 09] A. Nalli A, P. Haukkanen, \textit{On generalized Fibonacci and
Lucas polynomials}, Chaos, Solitons and Fractals, 42(5)(2009), 3179--86.%
\newline
[Ra; 15] J. L. Ramirez, \textit{Some Combinatorial Properties of the} $k$%
\textit{-Fibonacci and the} $k$\textit{-Lucas Quaternions}, An. St. Univ.
Ovidius Constanta, Mat. Ser., 23 (2)(2015), 201--212.\newline
[ Sa; 17] D. Savin, \textit{About Special Elements in Quaternion Algebras
Over Finite Fields}, Advances in Applied Clifford Algebras, June 2017, Vol.
27, Issue 2, 1801-1813.\newline
[St; 06] A.P. Stakhov, \textit{Fibonacci matrices, a generalization of the
\textquotedblleft Cassini formula\textquotedblright , and a new coding theory%
}, Chaos, Solitons and Fractals, 30(2006), 56-66.\newline
[Sw; 73] M. N. S. Swamy, \textit{On generalized Fibonacci Quaternions,} The
Fibonacci Quaterly, 11(5)(1973), 547--549.\newline

\begin{equation*}
\end{equation*}

\bigskip

Cristina FLAUT

{\small Faculty of Mathematics and Computer Science, Ovidius University,}

{\small Bd. Mamaia 124, 900527, CONSTANTA, ROMANIA}

{\small http://cristinaflaut.wikispaces.com/;
http://www.univ-ovidius.ro/math/}

{\small e-mail: cflaut@univ-ovidius.ro; cristina\_flaut@yahoo.com}

\medskip \medskip \qquad\ \qquad\ \ 

Diana SAVIN

{\small Faculty of Mathematics and Computer Science, }

{\small Ovidius University, }

{\small Bd. Mamaia 124, 900527, CONSTANTA, ROMANIA }

{\small http://www.univ-ovidius.ro/math/}

{\small e-mail: \ savin.diana@univ-ovidius.ro, \ dianet72@yahoo.com}\bigskip
\bigskip

\end{document}